\documentclass[10pt]{article}

\usepackage{amsmath, amssymb, amscd, fancyhdr}

\newcommand{\version}{version 1.0,\ \ Aug. 31, 2012}

\def\eqref#1{(\ref{#1})}
\newcommand{\goth}{\mathfrak}

\newcommand{\arrow}{{\:\longrightarrow\:}}

\newcommand{\C}{{\Bbb C}}
\newcommand{\R}{{\Bbb R}}

\def\1{\sqrt{-1}\:}
\newcommand{\restrict}[1]{{\left|_{{\phantom{|}\!\!}_{#1}}\right.}}
\newcommand{\cntrct}                
{\hspace{2pt}\raisebox{1pt}{\text{$\lrcorner$}}\hspace{2pt}}

\def\Bbb#1{\mathbb #1}

\renewcommand{\phi}{\varphi}
\renewcommand{\epsilon}{\varepsilon}

\newcommand{\sing}{\operatorname{\text{\sf sing}}}
\newcommand{\Vol}{\operatorname{Vol}}

\newcommand{\Tw}{\operatorname{Tw}}

\renewcommand{\Re}{\operatorname{Re}}
\renewcommand{\Im}{\operatorname{Im}}

\newcommand{\comment}[1]{{}}

\def\blacksquare{\hbox{\vrule width 4pt height 4pt depth 0pt}}
\def\endproof{\blacksquare}

\makeatletter

\newcounter{Mycounter}[section]
\newcounter{lemma}[section]
\newcounter{claim}[section]
\newcounter{sublemma}[section]
\newcounter{corollary}[section]
\renewcommand{\thecorollary}{\noindent{Corollary \thesection.\arabic{corollary}}}
\newcommand{\corollary}{%
     \setcounter{corollary}{\value{Mycounter}}
     \refstepcounter{corollary}
     \stepcounter{Mycounter}
     {\bf \thecorollary:\ }}

\newcounter{theorem}[section]
\renewcommand{\thetheorem}{\noindent{Theorem \thesection.\arabic{theorem}}}
\newcommand{\theorem}{%
     \setcounter{theorem}{\value{Mycounter}}
     \refstepcounter{theorem}
     \stepcounter{Mycounter}
     {\bf \thetheorem:\ }}

\newcounter{conjecture}[section]
\newcounter{proposition}[section]
\newcounter{definition}[section]
\renewcommand{\thedefinition}
       {\noindent{Definition~\thesection.\arabic{definition}}}
\newcommand{\definition}{%
     \setcounter{definition}{\value{Mycounter}}
     \refstepcounter{definition}
     \stepcounter{Mycounter}
     {\bf \thedefinition:\ }}

\newcounter{example}[section]
\newcounter{remark}[section]
\renewcommand{\theremark}{\noindent{Remark \thesection.\arabic{remark}}}
\newcommand{\remark}{%
     \setcounter{remark}{\value{Mycounter}}
     \refstepcounter{remark}
     \stepcounter{Mycounter}
     {\bf \theremark:\ }}

\newcounter{problem}[section]
\newcounter{question}[section]
\renewcommand{\thequestion}{\noindent{Question \thesection.\arabic{question}}}
\newcommand{\question}{%
     \setcounter{question}{\value{Mycounter}}
     \refstepcounter{question}
     \stepcounter{Mycounter}
     {\bf \thequestion:\ }}

\@addtoreset{equation}{section}
\@addtoreset{footnote}{section}
\makeatother

\begin{document}

\begin{center}
{\LARGE\bf Pseudoholomorphic curves on \\[3mm]
nearly K\"ahler manifolds}
\\[4mm]
Misha Verbitsky\footnote{Partially supported by RFBR grants
 12-01-00944-Á,  10-01-93113-NCNIL-a, and
AG Laboratory NRI-HSE, RF government grant, ag. 11.G34.31.0023.}
\\[4mm]

{\tt verbit@verbit.ru}
\end{center}

{\small 
\hspace{0.15\linewidth}
\begin{minipage}[t]{0.7\linewidth}
{\bf Abstract} \\
Let $M$ be an almost complex manifold equipped
with a Hermitian form such that its de Rham differential
has Hodge type (3,0)+(0,3), for example a nearly 
K\"ahler manifold. We prove that any connected component
of the moduli space
of pseudoholomorphic curves on $M$ is compact. This can be used
to study pseudoholomorphic curves on a 6-dimensional
sphere with the standard ($G_2$-invariant) almost 
complex structure.
\end{minipage}
}

{
\small
\tableofcontents
}


\section{Introduction}


Nearly K\"ahler manifolds 
were defined and studied by 
Alfred Gray, \cite{_Gray:structure_NK_}, in a general
context of intrinsic torsion of $U(n)$-\-struc\-tures
and weak ho\-lo\-no\-mies. An almost complex Hermitian
manifold $(M,I)$ is called {\bf nearly K\"ahler}, in this
sense, if $\nabla_X(I) X=0$, for any vector fields $X$ ($\nabla$ 
denotes the Levi-Civita connection). In other words, the
tensor $\nabla\omega$ must be totally skew-symmetric, for
$\omega$ the Hermitian form on $M$. If
$\nabla_X(\omega)\neq 0$ for any non-zero
vector field $X$, $M$ is called {\bf strictly nearly K\"ahler}.

For the last 10 years, 
the term ``nearly K\"ahler'' most often denotes
strictly nearly K\"ahler 6-manifolds. In sequel we shall follow
this usage, omitting ``strictly'' and ``6-dimensional''.

In dimension 6, a manifold is (strictly) nearly 
K\"ahler if and only if it admits a Killing spinor
(\cite{_Grunewald_}). Therefore, such a manifold is
Einstein, with positive Einstein constant. 

As one can easily show (see e.g. \cite{_Verbitsky:NK_}), 
strictly nearly K\"ahler 6-manifolds can be defined as
6-manifolds with structure group $SU(3)$ and fundamental
forms $\omega\in \Lambda^{1,1}_\R(M)$, 
$\Omega\in \Lambda^{3,0}(M)$, satisfying 
$d \omega=3\lambda\Re\Omega$,
$d \Im\Omega = -2 \lambda\omega^2$.
 An excellent
introduction to nearly K\"ahler geometry is
found in \cite{_Moroianu_Nagy_Semmelmann:NK_}.

Only 4 compact examples of nearly K\"ahler manifolds
are known, all of them homogeneous.
In \cite{_Butruille_} it was shown that any 
homogeneous nearly K\"ahler 6-manifold 
belongs to this list. Existence of non-homogeneous
compact examples is conjectured, but not proven yet. 

\begin{enumerate}

\item The 6-dimensional sphere $S^6$.  The 
almost complex structure on $S^6$ is reconstructed
from the octonion action, and the metric is standard.

\item $S^3\times S^3$, with the complex structure
mapping $\xi_i$ to $\xi'_i$,  $\xi_i'$ to $-\xi_i$,
where $\xi_i$, $\xi'_i$, $i=1,2,3$ is a basis
of left invariant 1-forms on the first and the second 
component.

\item Given a self-dual Einstein Riemannian 4-manifold $M$
with positive Einstein constant, one defines
its {\bf twistor space} $\Tw(M)$ as a total space of
a bundle of unit spheres in $\Lambda^2_-(M)$
of anti-self-dual 2-forms. Then $\Tw(M)$ 
has a natural K\"ahler-Einstein structure $(I_+, g)$,
obtained by interpreting unit vectors
in  $\Lambda^2_-(M)$ as  complex structure
operators on $TM$. Changing the sign of 
$I_+$ on $TM$, we obtain an almost complex structure $I_-$
which is also compatible with the metric $g$
(\cite{_Eels_Salamon_}). 
A straightforward computation insures that
$(\Tw(M), I_-, g)$ is nearly K\"ahler
(\cite{_Muskarov_}).

As N. Hitchin proved,
there are only two compact self-dual
Einstein 4-manifolds: $S^4$ and $\C P^2$.
The corresponding twistor spaces are
$\C P^3$ and the flag space $F(1,2)$.
The almost complex structure operator $I_-$ 
induces a nearly K\"ahler structure on these
two symmetric spaces. 

\end{enumerate}


\section{Pseudoholomorphic curves on nearly K\"ahler manifolds}


\subsection{Pseudoholomorphic curves on nearly K\"ahler manifolds}

The study of pseudoholomorphic curves on
almost complex manifolds is a big subject, spurred
by advances in physics and symplectic topology.
Even before the advent of the theory of 
pseudoholomorphic curves, complex curves in a 6-sphere
and a neary K\"ahler $\C P^3$ were used to
study the minimal surfaces. 

For pseudoholomorphic curves in a 6-sphere,
the basic reference is a paper
of R. Bryant \cite{_Bryant:octo_}.
Bryant observed that the complex curves
in $S^6$ can be approached similarly
to the real curves in $\R^3$, by computing
their first and second fundamental forms and the
torsion form. These operators, as defined by Bryant,
turn out to be holomorphic on a curve. Bryant
classifies the curves which have vanishing torsion
tensor, and constructs many examples of such curves,
using the correspondence between pseudoholomorphic
curves in $S^6$ and holomorphic curves in a 
5-dimensional quadric in $\C P^6$.

Results of Bryant were extended to $\C P^3$
with the nearly K\"ahler almost complex structure
by Feng Xu (\cite{_Feng_Xu:CP^3_}).
The usual (K\"ahler) projective 3-space $\C P^3$ 
can be obtained as a twistor space of $S^4$ fibered over $S^4$
with fibers isomorphic to $\C P^1$,$\C P^3\stackrel{\pi}\arrow S^4$ . This gives
a direct sum decomposition $T\C P^3=\pi^* T S^4 \oplus (\pi^* T S^4)^\bot$.
From the twistor construction, it is clear that this decomposition
is compatible with the complex structure.
To obtain the nearly K\"ahler complex structure,
we take the opposite complex structure on $(\pi^* T S^4)^\bot$,
and the usual one on $\pi^* T S^4$.

The pullback $\pi^* T S^4$ is 
identified with the holomorphic contact bundle on $\C P^3$
and the Legendrian curves (ones that are tangent to the
contact distribution) are by the above construction 
pseudoholomorphic with respect to the nearly K\"ahler 
almost complex structure.

Feng Xu shows that a pseudoholomorphic curve is Legendrian
if and only if its torsion vanishes, and classifies
pseudoholomorphic curves of genus 0.

\subsection{Compactness of the moduli spaces}

\definition
An {\bf almost complex Hermitian manifold}
is a manifold $(M,I,\omega)$ equipped with 
an almost complex structure $I$ and a
{\bf Hermitian form} $\omega$, that is, 
a (1,1)-form $\omega:=g(\cdot, I\cdot)$,
where $g$ is a {\bf Hermitian} ($I$-invariant
Riemannian) metric.

\hfill

\definition
Let $(M,I,\omega)$ be an almost complex Hermitian manifold.
A {\bf smooth, open pseudoholomorphic curve}
is a smooth 2-dimensional submanifold $S_0\subset M$
(not necessarily closed) satisfying $I(T_s S_0)=T_sS_0$
for each point $s\in S_0$.

\hfill

\remark
By Wirtinger's theorem, the Riemannian volume
form on a smooth pseudoholomorphic curve $S_0$
is equal to $\omega\restrict S_0$.

\hfill

\definition
Let $S\subset (M,I)$ be a compact subset,
smooth and 2-dimensional outside of a subset $S_{\sing}\subset S$
of Hausdorff dimension 0. Assume that 
the {\bf smooth part} $S \backslash S_{\sing}$
is pseudoholomorphic, and $S$ is its closure. 
Then $S$ is called {\bf a pseudoholomorphic
curve} in $M$. \footnote{Sometimes such curves are also
called {\bf holomorhic curves} and {\bf complex curves}.}

\hfill

\definition
Define {\bf the volume} of a pseudoholomorphic curve
as an integral $\int_{S\backslash S_{\sing}} \omega$.

\hfill

\remark
This integral equal to the Riemannian volume of $S\backslash S_{\sing}$.
Since $S_{\sing}$ has
Hausdorff codimension 2 in $S$, it is always finite (\cite{_Federer_}). 

\hfill

\definition
Let $M$ be a metric space, and ${\goth C}$ the set of its compact
subsets. Define {\bf the Hausdorff distance} on ${\goth C}$
as 
\[ 
  d_H(X,Y):=\max(\sup_{x\in X}(\inf_{Y\in Y}(d(x,y))),
  \sup_{y\in Y}(\inf_{x\in X}(d(x,y))).
\]
It is well known (see \cite{_Gromov:Riemannian_})
that $d_H$ defines a metric on ${\goth C}$,
which is complete and, when $M$ is a manifold, locally compact.

\hfill

The following fundamental theorem about pseudoholomorphic curves
is due to M. Gromov (\cite{_Gromov:curves_}, \cite{_Audin_}).

\hfill

\theorem\label{_Gromov_Compactness_Theorem_}
Let $(M,I,\omega)$ be a compact almost complex
Hermitian manifold, and ${\goth S}$ the set of
pseudoholomorphic curves on  $(M,I,\omega)$, equipped
with the Hausdorff distance.
Consider the volume function $\Vol:\; {\goth S} \arrow \R$.
Then $\Vol$ is continuous, and, moreover, $\Vol^{-1}([0,C])$
is compact, for any $C\in \R$.

\hfill

\remark
The set ${\goth S}$ with the Hausdorff distance and
the induced topology is called {\bf the moduli space 
of pseudoholomorphic curves}. It is a finite-dimensional,
locally compact topological space. The moduli of
pseudoholomorphic curves is a
real analytic variety if $(M,I,\omega)$ is
real analytic.

\hfill

The main result of this section is the following
theorem, which is then applied to the nearly K\"ahler
manifolds, such as $S^6$.

\hfill

\theorem\label{_volume_constant_Theorem_}
Let $(M,I,\omega)$ be an almost complex Hermitian
manifold satisfying $d\omega\in \Lambda^{3,0}(M)\oplus \Lambda^{0,3}(M)$.
Then the volume function $\Vol$ is constant on each connected
component of the moduli space of pseudoholomorphic curves.

\hfill

{\bf Proof:} Let $R$ be a connected 
component of the moduli space, and $\gamma:\; [0,1] \arrow R$
a continuous path. Denote by $R_\gamma\subset [0,1]\times M$
a set 
\[
R_\gamma:= \{(t,m)\in [0,1]\times m\ \ |\ \ m\in \gamma(t)\}.
\]
Without restricting the generality, we may assume that
$R_\gamma$ is smooth outside of a Hausdorff codimension 2
subset (see \cite{_Audin_}, where the local structure
of the moduli of pseudoholomorphic curves is describe 
analytically). Therefore, we may integrate bounded
differential forms over $R_\gamma$ and use the Stokes'
theorem as if $R_\gamma$ were a smooth, compact manifold
with boundary.

Denote by $\gamma:\; R_\gamma \arrow M$
the tautological projection. 
Then $\Vol(\gamma(1))-\Vol(\gamma(0))= \int_{\partial R_\gamma}\gamma^*\omega$.
Therefore, Stokes' theorem gives
\[
\Vol(\gamma(1))-\Vol(\gamma(0))=\int_{R_\gamma}\gamma^*d\omega.
\]
Consider the projection map $\pi:\; R_\gamma\arrow [0,1]$.
Then $\int_{R_\gamma}\gamma^*d\omega=\int_{[0,1]}\pi_*\gamma^*d\omega$,
where $\pi_*$ is a pushforward of differential forms. 

However, $\pi_*d\omega=0$, because $d\omega$ is a 
$(3,0)+(0,3)$-form, and $\pi_*$ is an integration over
$(1,1)$-cycles. \endproof

\hfill

Comparing this statement with Gromov's compactness
theorem (\ref{_Gromov_Compactness_Theorem_})
we obtain the following result (applicable to $S^6$ and other
nearly K\"ahler varieties).

\hfill

\corollary\label{_compactness_theorem_Corollary_}
Let $(M,I,\omega)$ be an almost complex Hermitian
manifold satisfying $d\omega\in \Lambda^{3,0}(M)\oplus \Lambda^{0,3}(M)$.
Then each connected component of the moduli of 
pseudoholomorphic curves on $M$ is compact.
Moreover, there is only a finite number of
components with volume bounded from above by any $C\in \R$.

\subsection{Open questions}

From \ref{_compactness_theorem_Corollary_}, it follows that
the pseudoholomorphic curves on nearly K\"ahler manifolds
behave in essentially the same way as on symplectic
manifolds: their moduli space is a countable union
of compact components. On K\"ahler manifolds the 
following questions are very easy to answer, but
on nearly K\"ahler manifolds, they seem quite mysterious.

\hfill

\question
With each component of the moduli of pseudoholomorphic
curves, one associates two numbers: a genus of a generic
curve in a family, and its volume, which is independent
from the choice of a curve by \ref{_volume_constant_Theorem_}.
Is there any relation between those numbers? What real 
numbers can occur as volumes of pseudoholomorphic
curves on a given nearly K\"ahler manifold?

\hfill

Notice that in \cite{_Bryant:octo_}, R. Bryant proved
existence of complex curves of arbitrary genus on $S^6$.

\hfill

It is interesting that on any nearly K\"ahler manifold $M$
the first Chern class $c_1(M)$ vanishes, but the curves behave
in a way which is completely different from observed on 
the Calabi-Yau manifolds.
In particular, the big-dimensional families of rational curves
seem to occur quite often, while on Calabi-Yau this is impossible
for topological reasons.

\hfill

\question
The moduli space of pseudoholomorphic curves on a
K\"ahler manifold is a complex manifold. What kind of geometric
structure arises on the  moduli space of pseudoholomorphic curves on a
nearly K\"ahler manifold?

\hfill

{\bf Acknowledgements:}
I am grateful to Klaus Hulek, Olaf Lechtenfeld and 
Alexander Popov for  discussions which lead to
the writing of this paper, and to Blaine Lawson for
an interesting discussion of this subject. Part of this
paper was written during my visit to the
Simons Center for Geometry and Physics; I wish
to thank SCGP for hospitality and stimulating atmosphere.

\hfill

{\small

\hfill

\noindent {\sc Misha Verbitsky\\
Laboratory of Algebraic Geometry, \\
Faculty of Mathematics, National Research University HSE,\\
7 Vavilova Str. Moscow, Russia
 }

}

\end{document}